\documentclass[11pt,reqno,oneside,final]{amsart}
\usepackage{graphicx, epsfig,cite,verbatim}
\usepackage{amsfonts,amsmath,amssymb}
\usepackage{float}
\usepackage{theorem}

\newcounter{remark} \def\openrem#1#2{\refstepcounter{remark}\bigskip
  {\noindent\sc#1~\theremark\if#2!{. }\else{ (#2). }\fi} }
\newenvironment{remark}[1][!]{\openrem{Remark}{#1}}{\endproof}

\def\R{\mathbb{R}}
\def\N{\mathbb{N}}
\def\Z{\mathbb{Z}}

\def\WW{{\rm W}}
\def\CC{{\rm C}}

\def\LL{{\rm L}}

\def\dx{\,{\rm d}x}

\def\pp{\partial}

\def\<{\langle}
\def\>{\rangle}

\def\eps{\varepsilon}

\def\E{\mathcal{E}}

\def\X{\mathcal{X}}

\def\Ls{\mathcal{L}}

\def\Cs{\mathcal{C}}

\renewcommand{\cases}[1]{\left\{ \begin{array}{rl} #1 \end{array} \right.}

\newcommand{\tworow}[2]{\genfrac{}{}{0pt}{}{#1}{#2}}

\newcommand{\smfrac}[2]{{\textstyle \frac{#1}{#2}}}
\def\half{\smfrac{1}{2}}
\def\quarter{\smfrac{1}{4}}

\begin{document}

\date{\today}

\title[Analysis of Node-Based Cluster Methods in the QC
  Method]{An Analysis of Node-Based Cluster Summation Rules in the Quasicontinuum
  Method}

\author{Mitchell Luskin}
\author{Christoph Ortner}

\address{Mitchell Luskin \\
School of Mathematics \\
University of Minnesota \\
206 Church Street SE \\
Minneapolis, MN 55455 \\
U.S.A.}
\email{luskin@umn.edu}

\address{Christoph Ortner\\
Oxford University
    Computing Laboratory\\
     Wolfson Building\\
     Parks Road, Oxford OX1
    3QD\\
     UK}
     \email {christoph.ortner@comlab.ox.ac.uk}

\thanks{
M. Luskin was supported in part by the NSF under
    grants DMS-0757355 and DMS-0811039, the DOE under Award Number
    DE-FG02-05ER25706, the Institute for Mathematics and Its
    Applications, and the University of Minnesota Supercomputing
    Institute.  C. Ortner was supported by EPSRC Critical Mass Project
    ``New Frontiers in the Mathematics of Solids''.
}



\begin{abstract}
  We investigate two examples of node-based cluster summation rules
  that have been proposed for the quasicontinuum method: a force-based
  approach (Knap \& Ortiz, J. Mech. Phys. Solids 49, 2001), and an
  energy-based approach which is a generalization of the non-local
  quasicontinuum method (Eidel \& Stukowski, J. Mech. Phys. Solids, to
  appear). We show that, even for the case of nearest neighbour
  interaction in a one-dimensional periodic chain, both of these
  approaches create large errors when used with graded and, more
  generally, non-smooth meshes. These errors cannot be removed by
  increasing the cluster size. We offer some suggestions how the
  accuracy of (cluster) summation rules may be improved.
\end{abstract}

\pagestyle{myheadings} \thispagestyle{plain} \markboth{MITCHELL LUSKIN
  AND CHRISTOPH ORTNER}{CLUSTER SUMMATION RULES IN THE QC METHOD}

\maketitle

\section{Introduction}
\label{sec:intro}
The quasicontinuum (QC) method
\cite{Ortiz:1995a,Miller:2003a,Shenoy:1999a,Shimokawa:2004} is a
prototypical coarse-graining technique for the static and quasi-static
simulation of crystalline solids. One of its key features is that,
instead of coupling an atomistic model to a continuum model, it uses
the atomistic model also in the continuum region where degrees of
freedom are removed from the model by means of piecewise linear
interpolation.

However, the nonlocal nature of the atomistic interactions makes
further approximation necessary to enable the computation of energies
or forces with complexity proportional to the number of coarse degrees
of freedom.  Two families of approximations have been developed to
achieve this goal.  One family of approximations localizes the
interactions by a strain energy density (based on the Cauchy--Born
rule) which provides sufficient accuracy in regions away from defects
where the strain gradient varies slowly.  Classical finite element
methodology can then be utilized in those regions modeled by the
strain energy density.  This class of quasicontinuum approximations
have been the subject of many recent mathematical analyses
~\cite{E:2006,LinP:2003a,LinP:2006a,E:2005a,Dobson:2008a,Dobson:2008c,Dobson:2008b,emingyang,Ortner:2008a}.

The purpose of the present paper is to investigate the second family
of approximations that have been developed to reduce the computational
complexity of the quasicontinuum method.  These methods, which have
received far less mathematical attention, use summation rules
(discrete variants of continuous quadrature rules) to approximate the
sums that define the quasicontinuum energy or forces.  To the best of
our knowledge, the force-based cluster summation rule of Knap and
Ortiz \cite{Knap:2001a}, the {\em non-local QC method} based on a
simple trapezoidal rule \cite[Sec. 3.3]{Miller:2003a}, or its
extension to energy-based cluster summation rules \cite{Eidel:2008a}
have not been analyzed to date.  These cluster summation rules
approximate the sum over atom-based quantities by uniformly averaging
over atoms in clusters (balls) around the nodes and then by weighting
these cluster averages so that summands that are obtained from
piecewise linear interpolation with respect to the quasicontinuum mesh
are exactly computed.

In a recent benchmark of different QC methods \cite{Miller:2008}, the
cluster summation rules do not compare favorably with quasicontinuum
approximations that utilize the strain energy density or with other
atomistic-to-continuum coupling methods. In the present paper, we give
a simple yet rigorous analytical explanation for this poor
performance. We demonstrate that, even for the simplest imaginable
atomistic model, a periodic one-dimensional chain with harmonic
nearest-neighbour interaction, the cluster summation rules formulated
in \cite{Knap:2001a,Eidel:2008a} lead to inconsistent and inaccurate QC
methods when used with graded and non-smooth meshes. Increasing the cluster size does not resolve this problem.
The benchmark \cite{Miller:2008} uses a mesh that is refined
from large triangles to the atomistic scale as is typical for
quasicontinuum computations.
Our analysis shows that this
kind of refined mesh would lead to an inconsistent and inaccurate method for
the approximation of our one-dimensional model by cluster summation rules.

The atomistic model, the finite element space ({\em coarse space}),
and some additional notation are introduced in Sections
\ref{sec:intro:model_problem} and \ref{sec:intro:cqc}. We treat the
two classes of cluster summation methods separately; the force-based
summation rule in Section \ref{sec:F} and the energy-based summation
rule in Section \ref{sec:E}. These sections can be read independently
of each other. Finally, in the conclusion, several possibilities are
identified how cluster summation rules might be modified in order to
lead to accuracte QC methods.

Although our analysis treats the approximation of the discrete sums in
the quasicontinuum energy, the reasons for the lack of accuracy of the
cluster summation rules can already be understood from applying the
cluster-method concepts to the finite element approximation of
continuum elasticity~\cite{Ciarlet:1978}.  The force conjugate to a
finite element nodal degree of freedom (the negative of the partial
derivative of the elastic energy functional constrained to the space
of finite element trial functions) depends only on integrals of the
jump of the displacement gradient across the element boundaries (the
finite element force density).  The accuracy of cluster-based
quadrature relies on the finite element force density being smooth in
space; however, the support of this force is concentrated on the
element boundaries. Thus, the {\em conjugate forces} obtained from the
cluster-based quadrature will be much too large.

Accurate node-based quadrature rules used in the approximation of a continuum
finite element energy are computed within each element during the assembly
process, which leads naturally to an accurate mesh-dependent non-uniform weighting
of the energy density in any ball surrounding each node.
The cluster-based quadrature approximation
uses a uniform weighting in the ball surrounding each node.  Thus,
since the energy density for a displacement in the finite element
trial space is generally discontinuous at the nodes, the cluster-based
quadrature rules are likely to be inaccurate for nonuniform meshes.

\begin{remark} {\normalfont We have not included the formulations of
    Lin \cite{LinP:2006a} or of Gunzburger \& Zhang
    \cite{Gunzburger:2008a,Gunzburger:2008b}, which are closely
    related to cluster summation methods, in our analysis. Our main
    reason for this exclusion is that their formulations do not suffer
    from the same deficiencies as the methods which we investiate
    here. {\em If} errors are present in their approach (we have not
    investigated this further), they would most likely be caused by
    finite range interaction and cannot be observed for the simple
    nearest-neighbour interaction system which we investigate here.}
\end{remark}

\subsection{The model problem}
\label{sec:intro:model_problem}
We choose the simplest immaginable atomistic model problem, a
one-dimensional periodic chain with nearest-neighbour pair potential
interaction. The continuum reference domain will be $(-1, 1]$. For
fixed $N \in \N$, the atomic spacing is given by $\eps = 1/N$, and the
atomistic reference lattice by
\begin{displaymath}
  \Ls = \big\{ \eps\ell : \ell = -N+1, \dots, N \big\}.
\end{displaymath}
The space of periodic displacements of $\Ls$ is denoted
\begin{equation}
  \label{eq:intro:space}
  \X = \big\{ v \in \R^\Z : v_{\ell+2N} = v_\ell \text{~for~} \ell \in \Z,
  \text{~and~} v_0 = 0 \big\}.
\end{equation}
We refer to Remark \ref{rem:intro:non-smooth_soln} for a motivation of
the constraint $v_0 = 0$.

We assume that the interatomic interaction reaches only nearest
neighbours, and that the only external force is a dead load. Thus, we
can write the {\em total energy} as a sum of a {\em stored energy}
$\E(v)$ and an {\em external potential energy} $-f[v]:$
\begin{align}
  \label{eq:intro:defn_energies}
  \Phi(v) =~&  \E(v) - f[v], 
 \intertext{where}
  \notag
  \E(v) =~& \sum_{\ell = -N+1}^N \eps\phi\big( \eps^{-1}(v_\ell
  - v_{\ell-1})\big), \quad \text{and} \\
  \notag
  f[v] =~& \sum_{-N+1}^N \eps f_\ell v_\ell.
\end{align}
Here $\phi$ is assumed to be smooth, at least in a neighbourhood of
zero, and $(f_\ell)_{\ell \in \Z}$ is a fixed $2N$-periodic
sequence. In fact, we shall usually make the simplifying assumption
that $\phi$ is a convex quadratic and that $f$ is obtained, for
example, by nodal interpolation from a smooth 2-periodic function
$\bar f(x)$. Such assumptions are valid for small displacements from
the reference state.

Rescaling the domain (and the energy) by the atomistic spacing $\eps$
is not strictly necessary, but it helps us understand the connection
of the atomistic problem to continuum theory.

The atomistic problem is to find
\begin{equation}
  \label{eq:intro:min_problem_a}
  u \in {\rm argmin}\,\Phi(\X),
\end{equation}
where `${\rm argmin}\,\Phi(\X)$' denotes the set of local minimizers
of $\Phi$ in $\X$. The first order necessary criticality condition, in
variational form, is
\begin{equation}
  \label{eq:intro:crit}
  \E'(u)[v] = f[v] \qquad \forall v \in \X,
\end{equation}
in short, $\E'(u) = f$.

\begin{remark}
  \label{rem:intro:non-smooth_soln}
  In the definition of the displacement space \eqref{eq:intro:space},
  we have imposed the condition $v_0 = 0$ for admissible
  displacements. This is one of several ways to remove the zero mode
  from the space in order to render the energy functional $\Phi$
  coercive.  Furthermore, this constraint allows us to easily
  construct a problem with a `singularity' at the origin (cf. Section
  \ref{sec:E:graded}). In general, if the external force $f$ is
  `smooth' and anti-symmetric, then the solution will be `smooth' as
  well. If $f$ is not anti-symmetric, then the solution may have a
  `kink' at the origin even if $f$ is smooth.
\end{remark}

We fix some additional notation, some of which we have already used
above. The arguments of nonlinear functionals are enclosed in round
brackets while those of (multi-)linear forms are enclosed in square
brackets, for example, $\E(u)$ or $f[u]$. The Fr\'{e}chet derivatives
are denoted by $'$, for example, $\E'(u)$ is a linear form on
$\X$. Consequently, $\E'(u)[v]$ denotes a directional
derivative. Similarly, $\E''(u)$ is a bilinear form on $\X$, and it is
written $\E''(u)[v, w]$ with arguments $v, w \in \X$. Finally, we will
frequently use the notation $v_\ell' = \eps^{-1}(v_\ell-v_{\ell-1})$
to denote the differences.

Atomistic displacements are always identified with their piecewise
affine interpolants. In particular, for $v \in \X$, we have $v_\ell' =
v'(x)$ for $x \in \left((\ell-1)\eps, \ell\eps\right)$. Through this
identification, the space $\X$ is naturally embedded in the spaces
\begin{displaymath}
  \WW^{1,p}_\#(-1, 1) = \big\{ v \in \WW^{1,p}(\R) : v(0) = 0,\
  v(x+2) = v(x) \quad \forall x \in \R \big\},
\end{displaymath}
for $1\le p \le \infty.$

\subsection{The constrained approximation}
\label{sec:intro:cqc}
The quasicontinuum approximation to the atomistic model problem
\eqref{eq:intro:min_problem_a} is obtained in two steps: (i) replacing
the displacement space $\X$ by a low-dimensional {\em coarse} space
$\X_h$, and (ii) approximating the nonlinear system
\eqref{eq:intro:crit} for arguments from the coarse space. Often, the
process is in fact reversed, however, for the class of QC methods
which we consider in the present paper, the order is as stated above.

We fix a set of {\em rep-atoms}
\begin{displaymath}
  \Ls_h = \big\{ \eps \ell_k : k = -K+1, \dots, K \big\} \subset \Ls,
\end{displaymath}
so that $\# \Ls_h \ll \# \Ls$. The set $\Ls_h$ is $2K$-periodically
extended, that is, we define $\ell_{k+2K}=\ell_{k}+2N$ for all $k \in
\Z$. The coarse space can therefore be written as
\begin{displaymath}
  \X_h = \big\{ v_h \in \X : v_h \text{~is~piecewise~affine~with respect to~}
  (\eps \ell_k)_{k \in \Z} \big\}.
\end{displaymath}

The {\em constrained atomistic approximation} is to find
\begin{equation}
  \label{eq:intro:CAA}
  \bar u_h \in {\rm argmin}\,\Phi(\X_h),
\end{equation}
for which the first order criticality condition is
\begin{equation}
  \label{eq:intro:CAA_crit}
  \E'(\bar u_h)[v_h] = f[v_h] \qquad
  \forall v_h \in \X_h.
\end{equation}
Even though the number of degrees of freedom is significantly
reduced in \eqref{eq:intro:CAA}, the nonlinear system
\eqref{eq:intro:CAA_crit} is still prohibitively expensive to evaluate
since it requires summation over {\em all} atoms. Hence, the second
step of the QC method, the approximation of the nonlinear system, is
as important as the coarsening step. One class of methods to achieve
this are the (cluster-)summation rules which we investigate in the
following sections \cite{Knap:2001a,Eidel:2008a}.

We conclude the introduction with some additional notation related to
the coarse space $\X_h$. For $k \in \Z,$ we denote $h_k = \eps(\ell_k
- \ell_{k-1})$. We will assume throughout that the mesh satisfies a
local regularity condition: there exists a constant $\kappa \geq 1$
such that
\begin{equation}
  \label{eq:intro:kmesh_condition}
  \kappa^{-1} h_{k-1} \leq h_k \leq \kappa h_{k-1} \qquad \text{for~}
  k = -K+1, \dots, K.
\end{equation}
For $v_h \in \X_h,$ we denote $V = (V_k)_{k \in \Z} = (v_{\ell_k})_{k
  \in \Z}$ the ($2K$-periodic) vector of nodal values, so that
\begin{equation}
  \label{eq:intro:basis_expansion}
  v_{h,\ell} = \sum_{k = -K+1}^K V_k \zeta_k(\eps\ell), \qquad
  \ell = -N+1, \dots, N,
\end{equation}
where $\zeta_k$ denotes the periodic nodal basis function associated
with node $\eps\ell_k$. Furthermore, we denote $V_k' = (V_k-V_{k-1}) /
h_k$ the {\em gradient} of $v_h$ in the element $(\eps\ell_{k-1}, \eps
\ell_k)$. In particular, we have $v_{h, \ell}' = V_k'$ if $\ell_{k-1}
< \ell \leq \ell_k$.

\section{Force-based summation rules}
\label{sec:F}
Since they are more easily understood, we shall first investigate
the force-based summation rules, introduced by Knap and Ortiz
\cite{Knap:2001a}. Our presentation closely follows their formulation.

Instead of viewing the constrained approximation \eqref{eq:intro:CAA}
as a minimization problem, we concentrate purely on the equilibrium
equations \eqref{eq:intro:CAA_crit}. However, instead of the
variational formulation \eqref{eq:intro:CAA_crit}, we use the {\em
  nodal force} formulation
\begin{equation}
  \label{eq:F:CAA_nodal}
  \frac{\pp \Phi(u_h)}{\pp U_j} = 0, \qquad j = -K+1, \dots, K.
\end{equation}
Using the expansion \eqref{eq:intro:basis_expansion} of $u_h$ in the
nodal basis $(\zeta_k)_{k = -K+1}^K$, the nodal forces are rewritten in
the form
\begin{displaymath}
  \frac{\pp \Phi(u_h)}{\pp U_j} = \sum_{\ell = -N+1}^N
  \frac{\pp \Phi(u)}{\pp u_\ell}\Big|_{u = u_h} \frac{\pp u_h(\eps\ell)}{\pp U_j}
  = \sum_{\ell = -N+1}^N \frac{\pp \Phi(u)}{\pp u_\ell}\Big|_{u = u_h} \zeta_j(\eps\ell),
\end{displaymath}
that is,
\begin{align}
  \label{eq:F:frx_expansion}
  \frac{\pp \Phi(u_h)}{\pp U_j} =~&
  \sum_{\ell = -N+1}^N F_\ell(u_h) \zeta_j(\eps\ell),
  \\ \intertext{where}
  \notag
  F_\ell(u) =~& \frac{\pp \Phi(u)}{\pp u_\ell}.
\end{align}
At this point, we apply a cluster summation rule to approximate the
sum in \eqref{eq:F:frx_expansion},
\begin{equation}
  \label{eq:F:FCQC}
  F_{j, h}(u_h) := \sum_{k = -K+1}^K \nu_k \sum_{\ell \in \Cs_k}
  F_\ell(u_h) \zeta_j(\eps \ell), \qquad j = -K+1, \dots, K,
\end{equation}
where the sets $\Cs_k$ are clusters surrounding the repatoms $\ell_k$,
\begin{displaymath}
  \Cs_k = \{ \ell_{k}-r_k^{-}, \dots, \ell_k + r_k^+ \},
  \qquad k = -K+1, \dots, K,
\end{displaymath}
and the weights $\nu_k$ are defined by the requirement that the basis
functions are summed exactly,
\begin{equation}
  \label{eq:F:weights}
  \sum_{\ell = -N+1}^N \zeta_j(\eps\ell) = \sum_{k = -K+1}^K \nu_k
  \sum_{\Cs_k} \zeta_j(\eps\ell), \qquad j = -K+1, \dots, K.
\end{equation}

In practice, the system \eqref{eq:F:weights} is solving using a {\em
  mass lumping} approximation \cite[Sec. 3.2]{Knap:2001a} yielding
approximate weights $\bar\nu_k$. We will only investigate the effect
of the cluster summation rule in the situation when $\eps \ll h_k$ for
all $k$, hence we shall assume throughout that $r_k^\pm \equiv r$ for
all $k$. In this particular case, we show in Appendix
\ref{sec:weights} that
\begin{displaymath}
  \bar\nu_k = \frac{h_k + h_{k+1}}{2 (2r+1) \eps},
  \qquad \text{and} \qquad
  \nu_k = \bar\nu_k + \mathcal{O}(1).
\end{displaymath}
We note that the {\em relative error} is of order
$\mathcal{O}(r\eps/(h_k + h_{k+1}))$.

\subsection{Analysis without external forces}
\label{sec:F:noforce}
In this section, we assume that $f \equiv 0$. To motivate this
assumption, we note that external body forces are only rarely the
driving force in an atomistic simulation and would therefore distort
the picture we are about to present. We shall simply ignore the fact
that, as a result, the atomistic problem becomes trivial.

If $f \equiv 0$, it follows that
\begin{displaymath}
  F_\ell(u) = \phi'(u_\ell') - \phi'(u_{\ell+1}'), \qquad
  \ell = -N+1, \dots, N,
\end{displaymath}
and, if we insert $u = u_h \in \X_h$, we obtain
\begin{equation}
  \label{eq:Fnn:frc_explicit}
  F_\ell(u_h) = \cases{ \phi'(U_k') - \phi'(U_{k+1}'), &
    \quad \text{if~} \ell = \ell_k, \\
    0, & \quad \text{otherwise.} }
\end{equation}
It follows that, independently of the cluster size, we obtain
\begin{displaymath}
  F_{k, h}(u_h) = \nu_k \big(\phi'(U_k') - \phi'(U_{k+1}')\big).
\end{displaymath}
Since $\nu_k \neq 0$, for all $k$, the equation $F_{k,h}(u_h) = 0$ is
equivalent to
\begin{displaymath}
  \frac{\pp\Phi(u_h)}{\pp U_k} = \phi'(U_k') - \phi'(U_{k+1}') = 0.
\end{displaymath}

Thus, we see that, even though the cluster summation rule
\eqref{eq:F:FCQC} is grossly inaccurate for moderate cluster radii
(the weights $\nu_k$ are of order $\mathcal{O}(h_k/(r\eps))$), the
resulting system is nevertheless equivalent to an exact evaluation of
the full constrained approximation which is known to be an excellent
approximation to the full atomistic system \cite{Ortner:2008a}.

\begin{remark}
  We need to be careful in extrapolating this observation to the case
  of finite range interaction and indeed the much more subtle and
  interesting two- and three-dimensional setting. These situations
  need to be investigated in more detail. Nevertheless, we can make
  some comments to motivate further investigation. The main
  observation that we have made in the present section, that forces
  are concentrated on the interfaces is still valid. In 2D and 3D, and
  for general atomistic models, the identity
  \begin{displaymath}
    \frac{\partial \Phi(u_h)}{\partial U_k} = \sum_{\ell \in \Ls}
    \frac{\pp \Phi(u)}{\pp u_\ell}|_{u = u_h} \zeta_k(\eps \ell)
  \end{displaymath}
  remains true (note, however, that now $\Ls \subset \R^d$, $d \in
  \{2, 3\}$, and $u : \Ls \rightarrow \R^d$). Furthermore, in the
  absence of an external force, in the interior of a large element,
  the force $\frac{\pp \Phi(u)}{\pp u_\ell}|_{u = u_h}$ is zero (or
  negligibly small) for an atomistic model with short range
  interaction. From this, we see that the contributions to the force
  on the node $\ell_k$ are concentrated {\em near} all element faces
  which touch the repatom $\ell_k$. This shows that the summation rule
  used to obtain $F_{k, h}(u_h)$ should be obtained from a summation
  over these faces (surface integration) instead of summation over the
  entire patch (volume integration).
\end{remark}

\subsection{Analysis with non-zero forces}
\label{sec:F:with_force}
We now return to the case of non-zero forces. We shall assume that the
forces $(f_\ell)_{\ell \in \Z}$ are obtained by interpolating a smooth
2-periodic function $\bar f \in \CC^2[-1,1]$. In this case, we obtain
\begin{displaymath}
  F_\ell(u) = \phi'(u_\ell') - \phi'(u_{\ell+1}') - \eps f_\ell, \qquad
  \ell = -N+1, \dots, N,
\end{displaymath}
and hence,
\begin{equation}
  \label{eq:Fnnn:frc_explicit}
  F_\ell(u_h) = \cases{ \phi'(U_k') - \phi'(U_{k+1}') - \eps f_{\ell_k}, &
    \quad \text{if~} \ell = \ell_k, \\
    - \,\eps f_\ell\,\;, & \quad \text{otherwise.} }
\end{equation}
It is then fairly straightforward to see that
\begin{displaymath}
  F_{k, h}(u_h) = \nu_k \big( \phi'(U_k') - \phi'(U_{k+1}') \big)
  - \tilde f_k,
\end{displaymath}
where $\tilde f_k$ is obtained by applying the cluster summation rule
to the external forces only,
\begin{displaymath}
  \tilde f_j = \sum_{k = -K+1}^K \nu_k \sum_{\ell \in \Cs_k} \eps f_\ell
  \zeta_j(\eps \ell), \qquad j = -K+1, \dots, K.
\end{displaymath}
Since we assumed that $\bar f \in \CC^2[-1,1]$, we can deduce from a
fairly straightforward interpolation error analysis (cf. Appendix
\ref{app:int_err}) that
\begin{equation}
  \label{eq:F:estimate_tildef}
  \tilde f_j = \sum_{\ell = -N+1}^N \eps f_\ell \zeta_j(\eps\ell)
  + \mathcal{O}(h_j^2 + h_{j+1}^2)
  = f[\zeta_j] + \mathcal{O}(h_j^2 + h_{j+1}^2).
\end{equation}
The force-based cluster quasicontinuum equations \eqref{eq:F:FCQC}
therefore become
\begin{displaymath}
  \nu_k \big( \phi'(U_k') - \phi'(U_{k+1}') \big) = f[\zeta_k]
  + \mathcal{O}(h_k^2 + h_{k+1}^2),
\end{displaymath}
as opposed to the `exact' equations of the constrained quasicontinuum
approximation
\begin{displaymath}
  \phi'(U_k') - \phi'(U_{k+1}') = f[\zeta_k].
\end{displaymath}

To illustrate this point further, let us assume that the interaction
is harmonic, that is $\phi(r) = \half r^2$, and that the mesh is uniform
($h_k = h$ for all $k$). In that case, the weights are given by $\nu_k
= h / (\eps (2r+1))$ (cf. Appendix \ref{sec:weights}), and hence
\begin{displaymath}
  \phi'(U_k') - \phi'(U_{k+1}') = \frac{\eps (2r+1)}{h} f[\zeta_k]
  + \mathcal{O}(h^2).
\end{displaymath}
Since the difference operator on the left-hand side is linear, we
therefore deduce that
\begin{displaymath}
  u_h = \frac{\eps (2r+1)}{h} \bar u_h + \mathcal{O}(h^2),
\end{displaymath}
where $\bar u_h$ is the solution of the constrained atomistic system
(\ref{eq:intro:CAA_crit}). In the typical case when 
$\eps r \ll h,$
this result demonstrates the catastrophic error made in the
force-based cluster summation rule. The reason why we do not observe a
similar cancellation effect as in Section \ref{sec:F:noforce} is
because the external force contribution was summed accurately, while
the summation of the interatomic forces is grossly inaccurate.

\section{Energy-based summation rules}
\label{sec:E}
We have seen in the previous section that the failure of the cluster
summation rule applied to the force balance equations fails because a
`volume integration' method was applied to a `surface integral'. It is
natural, therefore, to investigate the cluster summation rule applied
to the energy functional. This would lead to a conservative
coarse grained system, which was the main motivation for Eidel and
Stukowski \cite{Eidel:2008a} to use this method. They have noted in
\cite[Sec. 5]{Eidel:2008a} that this method also has shortcomings, and
we shall analyze these in detail in the present section.

To formulate the energy-based cluster summation rule, we first rewrite
the stored energy functional $\E$ in the form
\begin{align*}
  \E(u) =~& \sum_{\ell = -N+1}^N \eps \E_\ell(u), \quad
  \\ \intertext{where}
  \E_\ell(u) =~& \half \big(\phi(u_\ell') + \phi(u_{\ell+1}') \big).
\end{align*}
The term $\E_\ell(u)$ is the contribution of the atom at site $\ell$
to the overall energy. The sum over the terms $\E_\ell(u)$ is
approximated by a summation rule of the form
\begin{equation}
  \label{eq:E:csr_abstract}
  \sum_{\ell = -N+1}^N \eps g_\ell \approx \sum_{k = -K+1}^K \omega_k
  \sum_{\ell \in \Cs_k} g_\ell,
\end{equation}
where the sets $\Cs_k = \{ \ell_k - r_k^-, \dots, \ell_k + r_k^+\}$
are non-overlapping {\em clusters} surrounding the repatoms. The
weights $\omega_k$ are determined by requiring that the summation rule
is exact for all basis functions, that is,
\begin{equation}
  \label{eq:E:weight_rule}
  \sum_{\ell = -N+1}^N \eps \zeta_j(\eps\ell)
  = \sum_{k = -K+1}^K \omega_k \sum_{\ell \in \Cs_k} \zeta_j(\eps\ell),
  \qquad j = -K+1, \dots, K.
\end{equation}

To motivate a simplification which we are about to make, assume, for
the moment, that $r_k^\pm \equiv r$ for all $k$. For this case, we
have shown in Appendix \ref{sec:weights} that
\begin{equation}
  \label{eq:E:weights_explicit}
  \omega_k = \frac{h_k + h_{k+1}}{2 (2r+1)} + \mathcal{O}(\eps).
\end{equation}
Furthermore, we observe that
\begin{align*}
  \sum_{\ell \in \Cs_k} \E_\ell(v_h) =~& r \phi(V_k')
  + \half \big( \phi(V_k') + \phi(V_{k+1}') \big)
  + r \phi(V_{k+1}') \\[-2mm]
  =~& \half(2r+1) \big( \phi(V_{k}') + \phi(V_{k+1}') \big) \\[2mm]
  =~& (2r+1) \E_{\ell_k}(v_h).
\end{align*}
Thus, we see that, up to an error of order $\mathcal{O}(\eps)$, a
finite cluster size reduces immediately to a discrete trapezoidal
rule.

In view of this observation, we shall assume throughout this section
that $\Cs_k = \{\ell_k\}$. The approximate energy functional becomes
\begin{equation}
  \label{eq:E:qc_energy}
  \E_{h}(v_h) = \sum_{k = -K+1}^K \omega_k \E_{\ell_k}(v_h),
\end{equation}
with weights $\omega_k = \half(h_k + h_{k+1})$. This method
\eqref{eq:E:qc_energy} is sometimes labelled the {\em non-local QC
  method} \cite[Sec. 3.3]{Miller:2003a}.

\begin{remark}
  1. The observations made above are only partially valid for
  non-nearest neighbour interaction. In that case, additional
  interface terms of the form $\phi(V_k' + V_{k+1}')$ enter the QC
  energy functional.

  2. A further correction from our simplifying assumption needs to be
  taken into account when the mesh is refined to atomistic level where
  we need to use variable cluster sizes. For simplicity, we have
  chosen to ignore this further complication, but note that our
  analysis in Appendix \ref{sec:weights} can be generalized to
  variable cluster sizes provided the cluster radii are symmetric in
  each element (that is, $r_{k-1}^+ = r_k^-$). In that case, we would
  obtain a similar formula as \eqref{eq:E:weights_explicit} but with
  an error of order $\mathcal{O}(\eps \max_k r_k^\pm)$.
\end{remark}

\medskip For simplicity, we assume that the dead load $f[v]$ is not
approximated. The total energy for the QC method is therefore given by
\begin{displaymath}
  \Phi_{h}(v_h) = \E_{h}(v_h) - f[v_h],
\end{displaymath}
where $\E_h$ is defined in \eqref{eq:E:qc_energy}. The corresponding
criticality condition is
\begin{equation}
  \label{eq:E:crit_qc}
  \E_h'(u_h)[v_h] = f[v_h] \qquad \forall v_h \in \X_h.
\end{equation}

In order to analyze the (non-local) QC method, we first rewrite $\E_h$
in the form
\begin{displaymath}
  \E_h(v_h) = \sum_{k = -K+1}^K \omega_k \half \big( \phi(V_k')
  + \phi(V_{k+1}') \big)
  = \sum_{k = -K+1}^K \half(\omega_k + \omega_{k-1}) \phi(V_k').
\end{displaymath}
Using $\omega_k = \half(h_k + h_{k+1})$, we see that
\begin{align*}
  \half(\omega_k + \omega_{k-1}) =~& \quarter( h_{k-1}+h_k + k_k + h_{k+1} ) \\
  =~& h_k + \quarter( h_{k-1} - 2 h_k + h_{k+1} ) \\
  =:~& h_k (1 + \hat \omega_k),
\end{align*}
where
\begin{equation}
  \label{eq:3}
  \hat\omega_k = \frac{h_{k-1} - 2 h_k + h_{k+1}}{4 h_k},
\end{equation}
and hence, we obtain
\begin{equation}
  \label{eq:E:Eqc_error}
  \E_h(v_h) = \E(v_h) + \sum_{k = 1}^K h_k \hat\omega_k \phi(V_k').
\end{equation}
We use $\hat\omega$ to denote the piecewise constant function taking
values $\hat\omega_k$ in the elements $(\eps\ell_{k-1}, \eps \ell_k)$.

We can relate the connection between the error in the cluster
approximation to the error in the trapezoidal rule by noting that
\eqref{eq:E:Eqc_error} is equivalent to
\begin{equation}
  \label{eq:E:Eqc_error2}
  \E_h(v_h) = \E(v_h) + \quarter \sum_{k = 1}^K h_k
  \left(\phi(V_{k+1}')-2\phi(V_k')+\phi(V_{k-1}')\right).
\end{equation}

From \eqref{eq:E:Eqc_error} and \eqref{eq:E:Eqc_error2}, we already
anticipate that the local mesh smoothness will have a significant
impact on the accuracy of the method. For example, if the
$\hat\omega_k$ are oscillatory, then it is possible to lower the
energy by introducing a `microstructure' into the quasi-continuum
displacement. We will see this in detail in Example 3 below. In
Example 2, we discuss another effect that may introduce large errors
in the simulation.

\begin{remark}
  We could analyze the consistency of the method using finite
  difference techniques. Taking the derivative of $\E_h(u_h)$ with
  respect to the nodal value $U_k$, we obtain
  \begin{displaymath}
    \frac{\pp \E_h(u_h)}{\pp U_k} =
    \big(1 + \hat\omega_k\big) \phi'(U_k')
    - \big(1 + \hat\omega_{k+1}\big) \phi'(U_{k+1}')
    = f[\zeta_k],
  \end{displaymath}
  where $\hat\omega_k$ is given by \eqref{eq:3}. One can see here as
  well that, if the terms $\hat\omega_k$ are not close to zero, then
  the method is inconsistent. However, a rigorous error analysis is
  more conveniently performed in the variational setting of the finite
  element method.
\end{remark}

\medskip To further simplify the analysis, we assume from now on that
the interaction is harmonic, that is, $\phi(r) = \half r^2$. This
assumption can be justified, for example, for small perturbations from
a reference state. In that case, the fully atomistic problem
\eqref{eq:intro:crit} has a unique solution $u$. Furthermore, let
$\bar u_h$ be the unique solution of the constrained approximation,
which is the {\em best approximation} to $u$ from $\X_h$ in the energy
norm. Since all weights $\omega_k$, $k = 1, \dots, K,$ are positive,
it follows that the QC functional $\Phi_h$ also has a unique critical
point, $u_h$.

Since $\E$ and $\E_h$ are both quadratic, their Hessians are
independent of the arguments. Thus, we will write $\E''_{(h)}[v_h,
w_h]$ instead of, say, $\E''(u_h)[v_h, w_h]$. With this notation, the
criticality conditions \eqref{eq:intro:CAA_crit} and
\eqref{eq:E:crit_qc} become, respectively,
\begin{align*}
  \E''[\bar u_h, v_h] =~& f[v_h] \qquad \forall v_h \in \X_h, \quad
  \\
  \E_h''[u_h, v_h] =~& f[v_h] \qquad \forall v_h \in \X_h.
\end{align*}
Thus, the error $u_h - \bar u_h$ satisfies, for all $v_h
\in \X_h$,
\begin{displaymath}
  \E_h''[u_h - \bar u_h, v_h] = f[v_h] - \E_h''[\bar u_h, v_h]
  = \E''[\bar u_h, v_h] - \E_h''[\bar u_h, v_h].
\end{displaymath}
Using Strang's First Lemma \cite[Thm. 4.1.1]{Ciarlet:1978} and the
mesh regularity condition \eqref{eq:intro:kmesh_condition}, we obtain
\begin{align}
  \label{eq:E:strang}
  \half(1+\kappa^{-1}) \|u_h' - \bar u_h' \|_{\LL^2}
  \leq~& \big(\E_h''[u_h-\bar u_h, u_h - \bar u_h]\big)^{1/2}  \\
  \notag
  \leq~& \sup_{v_h \in \X, \|v_h'\|_{\LL^2} = 1} \big| \E''[\bar u_h,v_h]
  - \E_h''[\bar u_h, v_h] \big|.
\end{align}
We wish to obtain sharp upper and lower bounds on the {\em variational
  crime}. To this end, we first note that a piecewise constant
function $r_h$, is the gradient of an element of $\X_h$ if, and only
if, $\int_{-1}^1 r_h \dx = 0$.  Thus, setting \begin{equation}
  \label{eq:E:defn_a}
  a = \half \sum_{k = -K+1}^K h_k \hat\omega_k \bar U_k'
  = \half \int_{-1}^1 \hat\omega \bar u_h' \dx,
\end{equation}
we obtain
\begin{align}
  \notag
  \sup_{\tworow{v_h \in \X_h}{\|v_h'\|_{\LL^2} = 1}}
  \big| \E''(\bar u_h, v_h) - \E_h''(\bar u_h, v_h) \big|
  =~& \sup_{\tworow{v_h \in \X_h}{\|v_h'\|_{\LL^2} = 1}}
  \Big| \int_{-1}^1 \big( \hat\omega \bar u_h' - a\big) v_h' \dx \Big| \\
  \label{eq:E:defn_rho}
  =~& \| \hat \omega \bar u_h' - a\|_{\LL^2} =: \rho(\bar u_h),
\end{align}
which gives an upper bound on the error. To obtain a lower bound as well, we
reverse the argument in \eqref{eq:E:strang}, yielding
\begin{align*}
  \rho(\bar u_h) =~& \sup_{v_h \in \X_h, \|v_h'\|_{\LL^2} = 1}
  \big| \E''[\bar u_h, v_h]
  - \E_h''[\bar u_h, v_h] \big|  \\
  =~& \sup_{v_h \in \X_h, \|v_h'\|_{\LL^2} = 1}
  \big| \E_h''[u_h - \bar u_h, v_h] \big| \\
  \leq~& \half (1+\kappa) \|\bar u_h' - u_h' \|_{\LL^2}.
\end{align*}
Combining this estimate with \eqref{eq:E:strang} and \eqref{eq:E:defn_rho},
we finally arrive at
\begin{align}
  \label{eq:E:qc_error}
  &\half (1 + \kappa^{-1}) \| \bar u_h' - u_h' \|_{\LL^2} \leq \rho(\bar u_h)
  \leq \half(1+ \kappa) \| \bar u_h' - u_h' \|_{\LL^2}, \\
   \intertext{where} &\quad
  \rho(\bar u_h)^2 = \sum_{k = -K+1}^K h_k
  \big|\hat\omega_k \bar U_k' - a \big|^2
  = \| \hat\omega \bar u_h' - a \|_{\LL^2}^2\notag
\end{align}
Note that \eqref{eq:E:qc_error} does not estimate the actual error $u
- u_h$ but the deviation from the best approximation in the energy
norm. In the following, we will investigate the term $\rho(\bar u_h)$ for
three typical meshes that may occur in a simulation.

\subsection{Example 1: smooth meshes} We begin by looking at a
somewhat idealistic situation. We assume that $\eps \ll h$ and that
the mesh nodes at $\eps \ell_k$, $k = -K+1, \dots, K,$ are given by a
smooth (and periodic) deformation $\varphi$ of the periodic domain
$(-1, 1]$, that is,
\begin{displaymath}
  \eps \ell_k = \varphi(h k), \qquad k = -K+1, \dots, K,
\end{displaymath}
where $h=1/K$, and $\varphi(x) = \varphi(x)+2$ for all $x \in \R$. In
that case, the term $\hat\omega_k$ can be estimated, using Taylor's
Theorem, to obtain
\begin{align*}
  \hat\omega_k =~& \frac{h_{k+1} - 2 h_k + h_{k-1}}{4 h_k} \\
  =~& \frac{ \varphi(h(k+1)) - 3 \varphi( hk)
  + 3 \varphi(h(k-1)) - \varphi(h(k-2))}{ 4 (\varphi(hk) - \varphi((h-1)k))} \\
  =~& \quarter h^2 \frac{\varphi'''(\bar x_k)}{\varphi'(\bar x_k)}
  + \mathcal{O}(h^3),
\end{align*}
where $\bar x_k = (k-\half) h$. In particular, we obtain
\begin{displaymath}
  |\hat\omega_k| \leq C h^2 + \mathcal{O}(h^3),
\end{displaymath}
where $C = \quarter \max_{x \in [-1,1]} |\varphi'''(x)/\varphi'(x)|$.
Since
\begin{displaymath}
  \rho(\bar u_h)^2 =\sum_{k = -K+1}^K h_k \big|
  \hat\omega_k \bar U_k' \big|^2 - 2 a^2,
\end{displaymath}
it follows that
\begin{displaymath}
  \rho(\bar u_h)^2 \leq \sum_{k = -K+1}^K h_k \big|
  \hat\omega_k \bar U_k' \big|^2
  \leq C^2 h^4 \| \bar u_h' \|_{\LL^2}^2 + \mathcal{O}(h^6),
\end{displaymath}
which is of a smaller order than the best approximation error.

\subsection{Example 2: graded meshes}
\label{sec:E:graded}
Since the main target of the quasicontinuum method are problems with
defects or singularities, extremely smooth meshes satisfying the
assumptions of Example 1 are rare in quasicontinuum
applications. The most important example for the QC method is
a mesh which refines to atomistic level. To investigate this
situation, we construct an exponentially graded mesh as follows. We
fix $K > 0$ and $N = 2^{K-1}$, and define $\ell_0 = 0$ and
\begin{displaymath}
  \ell_k = {\rm sgn}(k) 2^{|k|-1}, \quad k = -K+1, \dots, K.
\end{displaymath}
In that case, we obtain
\begin{displaymath}
  h_k = \cases{
    \eps, & k = 0, 1, \\
    2^{k-2} \eps, & k = 2, \dots, K, \\
    2^{|k|-1} \eps, & k = -1, -2, \dots, -K+1, }
\end{displaymath}
which, in particular, gives the mesh regularity parameter $\kappa =
2$. We can, furthermore, explicitly compute the coefficients
\begin{displaymath}
  \hat\omega_k =  \cases{
  0, & k = 0, 1, \\
  1/4, & k = -1, 2, \\
  1/8, & k = -K+2,  \dots, -2, 3, \dots, K-1, \\
  -1/8, & k = -K+1, K. }
\end{displaymath}

To further investigate the error, let us assume that the displacement
gradient in the `continuum region' is negligable. Let us further
assume that the displacement gradient does not vary considerably
between the elements $(0, h)$ and $(h, 2h)$ as well as between $(-h,
0)$ and $(-2h, -h)$. In that case, we can ignore the $\hat\omega_K =
\hat\omega_{-K+1} = -\frac{1}{8}$ coefficients in the outmost
elements, and we can `split' the coefficients $\hat\omega_2 =
\hat\omega_{-1}$ among the purely atomistic elements in order to
obtain $a \approx 0$ and
\begin{displaymath}
  \rho(\bar u_h)^2 = \sum_{k = -K+1}^K h_k \big| \hat\omega_k \bar U_k'-a\big|^2
  \approx \smfrac{1}{8^2} \sum_{-K+1}^K h_k \big|\bar U_k'\big|^2
  = \big(\smfrac{1}{8} \|\bar u_h'\|_{\LL^2}\big)^2.
\end{displaymath}

\begin{figure}
  \begin{center}
    \includegraphics[width=8cm]{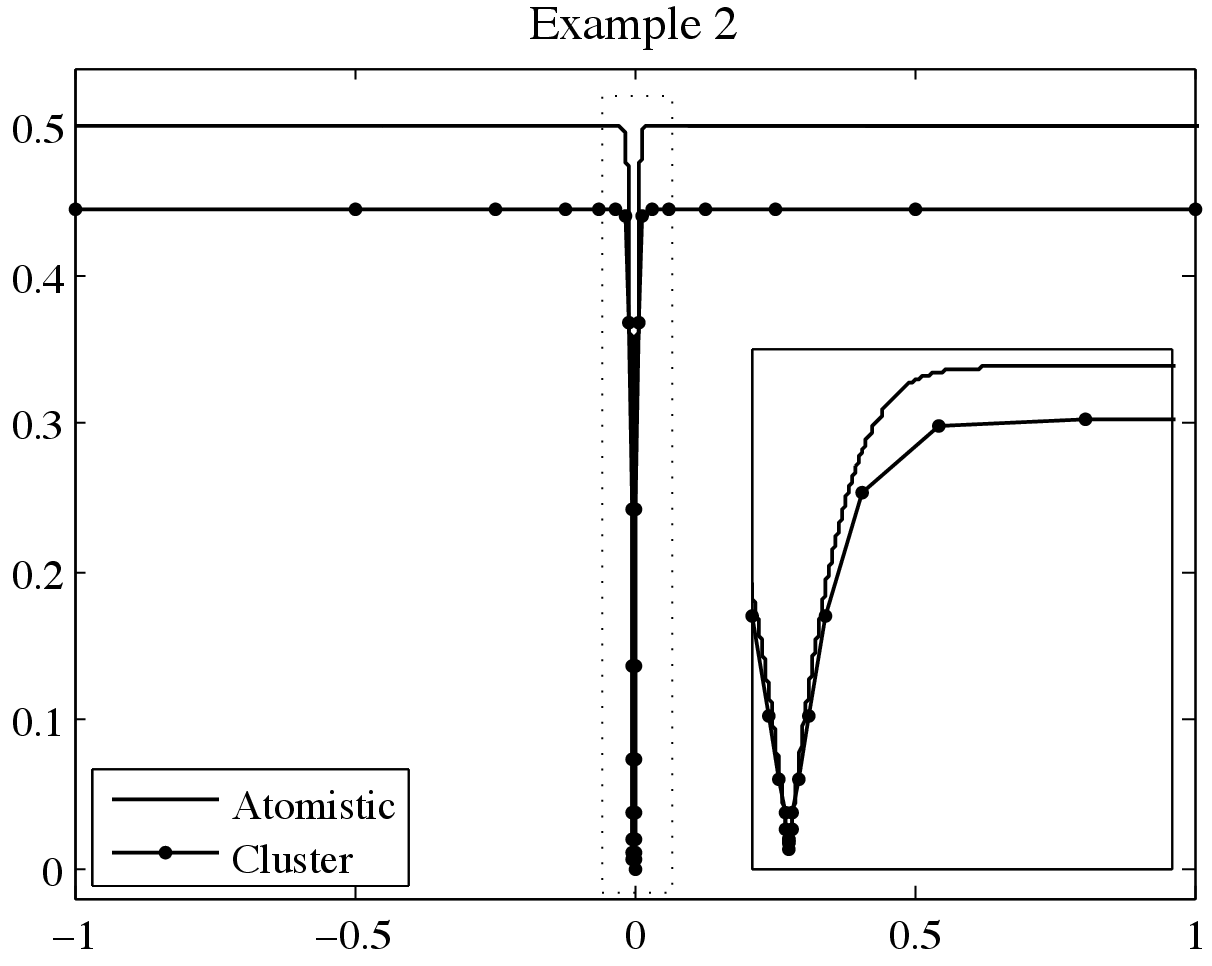}
  \end{center}
  \caption{\label{fig:example2} Computational example on a highly
    graded mesh with force $f(x) = 10^4 \exp(-10^4 x^2)$, $N =
    2^{14},$ and $K = 15$. Since $f$ is not anti-symmetric, the
    solution is non-smooth at the origin (cf. Remark
    \ref{rem:intro:non-smooth_soln}). The relative error in the energy
    norm satisfies $\| u_h' - \bar u_h' \|_{\LL^2}/\|\bar
    u_h'\|_{\LL^2} \approx 0.11$ (that is $11\%$), which is in
    excellent agreement with our prediction. The relative error for
    the energy satisfies $(\E(u) - \E_h(u_h)) / |\E(u)| \approx -0.13$
    (that is $13\%$). Precisely as we predicted, we see that the
    failure in the coefficients enforces a smaller QC displacement
    which results in a higher energy.}
\end{figure}

From \eqref{eq:E:qc_error}, we would therefore expect a relative error
of size $1/8$ (that is $12\%$), independent of the mesh size $h$. This
is in perfect agreement with the computational example we present in
Figure \ref{fig:example2}.

\subsection{Example 3: a non-smooth mesh}
In the final example of our analysis of the energy-based cluster
summation rule, we consider a mesh which is quasi-uniform, but {\em
  not} smooth. We assume that $\eps \ll h$, and that
\begin{equation}
  \label{eq:E:osc_mesh}
  h_{k} = \half h \big( 3 + (-1)^k \big), \qquad k = -K+1, \dots, K,
\end{equation}
that is, we have $h_1 = h$, $h_2 = 2h$, $h_3 = h$, and so forth. A
mesh of precisely this type will rarely be found in practise, however,
it is an excellent model situation that demonstrates a source of error
for non-smooth meshes.

In this situation, the coeffcients $\hat\omega_k$ satisfy
\begin{displaymath}
  \hat\omega_k = \cases{
    -1/4, & \text{if~} k \text{~is~even,} \\
    1/2, & \text{if~} k \text{~is~odd.}
    }
\end{displaymath}
Suppose that $\bar u_h$ is the interpolant of a smooth function $\bar
u$, so that $\bar U_k'$ varies little between elements. Then the
oscillatory nature of the coefficients $\hat\omega_k$, weighted
according to the size of the elements, indicates that $a \approx
0$. More precisely, we have
\begin{displaymath}
  h_k \hat\omega_k \bar U_k' = \cases{
    \smfrac{1}{6}(h_k \bar U_k' + h_{k+1} \bar U_{k+1}')
    + \smfrac{1}{6} h_{k+1} (\bar U_{k+1}' - \bar U_{k}'),
    & k \text{~odd}, \\
    -\smfrac{1}{6}(h_{k-1} \bar U_{k-1}' + h_k \bar U_k')
    + \smfrac{1}{6} h_k (\bar U_{k}' - \bar U_{k-1}'),
    & k \text{~even},
  }
\end{displaymath}
from which we can easily deduce that $|a| \leq C h$, where $C$ depends
on the second differences of $\bar U_k$ (which we assumed to be
moderately small). Some similar, algebraic manipulations show that
\begin{displaymath}
  \sum_{k = -K+1}^K h_k \big|\hat\omega_k \bar U_k'\big|^2 =
  \smfrac{1}{8} \| \bar u_h' \|_{\LL^2}^2
  + \mathcal{O}(h),
\end{displaymath}
and thus, we obtain
\begin{align*}
  \rho(\bar u_h)^2 = \sum_{k = -K+1}^K h_k \big|\hat\omega_k \bar U_k'\big|^2
  - 2a^2 \geq \smfrac{1}{8} \| \bar u_h' \|_{\LL^2}^2
  - \mathcal{O}(h).
\end{align*}

\begin{figure}
  \begin{center}
    \includegraphics[width=8cm]{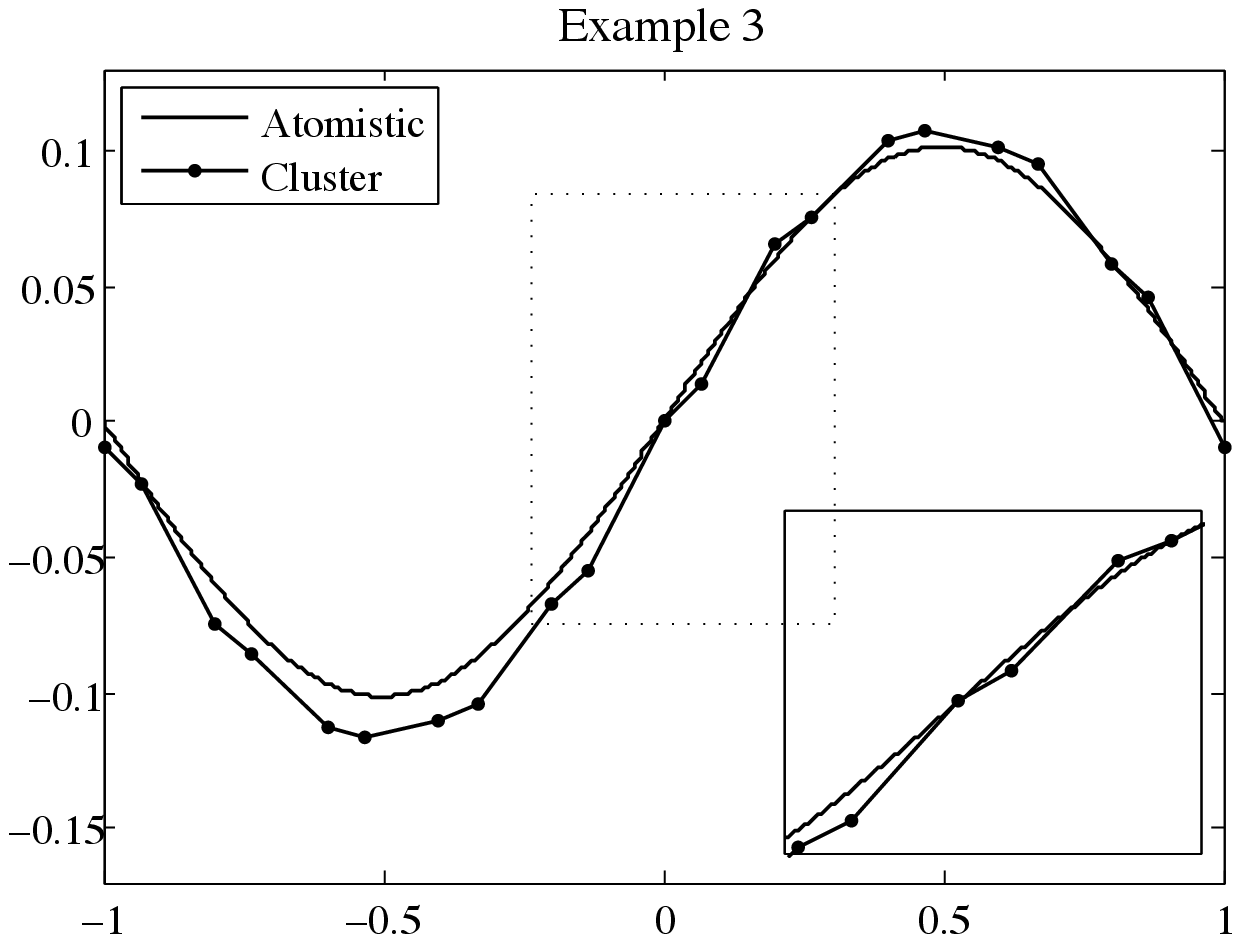}
  \end{center}
  \caption{\label{fig:example3} Computational example on an
    oscillatory mesh with force $f(x) = \sin(\pi x)$, $N = 10^4$ and
    $K = 20$. The fully atomistic solution is given by $u(x) =
    \pi^{-2} \sin(\pi x)$. The relative error in the energy norm
    satisfies $\| u_h' - \bar u_h' \|_{\LL^2}/\|\bar u_h'\|_{\LL^2}
    \approx 0.33$ (that is $33\%$) while the relative error for the
    energy satisfies $(\E(u) - \E_h(u_h)) / |\E(u)| \approx 0.097$
    (that is $9.7\%$). In the zoomed box, we see that the
    microstructure, induced by the oscillatory coefficients, is
    lowering the energy and creates a `non-smooth' quasicontinuum
    solution. }
\end{figure}

From our general error estimate \eqref{eq:E:qc_error}, we therefore
expect the relative error to be roughly of the order $1/\sqrt{8}$
(that is $35\%$), which is in excellent agreement with the numerical
example shown in Figure \ref{fig:example3}.

\section{Conclusion}
We have shown that node-based cluster summation rules, applied either
to the force-based formulation of the QC method or to the energy-based
formulation of the QC method lead to  inconsistent and inaccurate numerical
schemes when used with graded or non-smooth meshes. We stress, furthermore, that increasing the cluster size is
{\em not} a remedy for the sources of error which we have discussed.

We do not rule out, however, that QC methods based on a more careful
choice of summation points may yet lead to excellent computational
tools. We would like to comment on three options which qualify for
further investigation:

Lin's formulation \cite{LinP:2006a} and the formulation of Gunzburger
\& Zhang \cite{Gunzburger:2008a,Gunzburger:2008b}, which are based on
summation points in the interior of the elements do not suffer from
any of the deficiencies which we have found in the present work. It
will be necessary, however, to carefully investigate the effect of
next-nearest neighbour and finite range interaction in the transition
region in which the mesh is refined from large triangles to the
atomistic scale where all degrees of freedom are retained.

A force-based formulation, where the summation is performed over
element interfaces rather than elements may yet lead to an accurate QC
method. This is clearly true in one dimension, but needs to be
carefully studied in higher dimensions.

Finally, we propose to investigate the possibility of assigning
variable weights to atoms within the same cluster, in an energy-based
cluster summation rule. It can be readily verified, following the
analogy with continuum finite element energies discussed at the
end of the introduction, that if the
average over atoms in a cluster is weighted according to element
sizes, then the resulting method will be accurate for
nearest-neighbour interaction. Once again, the crucial questions are
whether this accuracy can be retained for finite range interaction and
the application relevant two- three-dimensional situations.

\appendix

\section{Proofs}

\subsection{Computation of summation weights}
\label{sec:weights}
The analysis presented in this appendix applies to both the
energy-based and the force-based summation rules, since the weights
satisfy $\omega_j = \eps \nu_j$. For no particular reason, we chose to
work with the weights $(\omega_j)_{j = -K+1}^K$.

We assume throughout that $r_k^\pm \equiv r$. According to the
requirement \eqref{eq:E:weight_rule}, the governing equations for the
weights $\omega = (\omega_k)_{k = -K+1}^K$ are $M\omega = g$ where
\begin{align*}
  (M \omega)_j :=~& \sum_{k = -K+1}^K \omega_k \sum_{\ell \in \mathcal{C}_k} \zeta_j(\eps\ell) \\
  =~& \omega_{j-1} \big( \half r(r+1) \eps h_j^{-1} \big)  + \omega_{j+1} \big( \half r(r+1) \eps h_{j+1}^{-1} \big) \\
  & + \omega_j \big( (2r+1) - \half r(r+1) \eps h_j^{-1} - \half r(r+1) \eps h_{j+1}^{-1}\big)
\end{align*}
and
\begin{displaymath}
  g_j := \sum_{\ell = -N+1}^N \eps \zeta_j(\eps\ell) = \half( h_{j+1} + h_j).
\end{displaymath}
To prove that $M$ is invertible, we show that it is row-diagonally
dominant. For each $j,$ we have
\begin{displaymath}
   M_{jj} - \sum_{k \neq j} |M_{j,k}|
  = (2r + 1) - r (r+1) \eps (h_j^{-1} + h_{j+1}^{-1}).
\end{displaymath}
Since we assumed that the clusters do not overlap, it follows that $\ell_j
- \ell_{j-1} \geq 2r+1$, in particular, $- \eps r > -\half h_j$, from which
we deduce that
\begin{equation}
  \label{eq:app:rdd}
  M_{jj} - \sum_{k \neq j} |M_{j,k}| > (2r+1) - (r+1) = r.
\end{equation}
Thus, $M$ is invertible and $\omega$ is well-defined. We note,
furthermore, that \eqref{eq:app:rdd} implies that
\begin{equation}
  \label{eq:app:normMinv}
  \| M^{-1} \|_{\ell^\infty} \leq \cases{ 1/r, & \text{if~} r \geq 1, \\
    1, & \text{if~} r = 0. }
\end{equation}

Our next observation is that the {\em lumped system} for computing the
approximate weights $\bar\omega = (\bar \omega_k)_{k=-K+1}^K$ is
\begin{displaymath}
  (2r+1) \bar \omega_j = \half (h_{j+1} + h_j), \qquad j = -K+1, \dots, K,
\end{displaymath}
that is,
\begin{equation}
  \label{eq:app:lumped_weights}
  \bar\omega_j = \frac{h_{j+1} + h_j}{2 (2r+1)}, \qquad j = -K+1, \dots, K.
\end{equation}

We shall now prove that the exact weights $(\omega_k)_{k = -K+1}^K$
are only $\mathcal{O}(\eps)$ perturbations from the approximate
weights obtained by mass-lumping. To this end, we define the residual
$\rho = (\rho_k)_{k = -K+1}^K$,
\begin{align*}
  \rho_j :=~& (M\bar\omega)_j - g_j \\
  =~& (\bar\omega_{j-1}-\bar\omega_{j}) \big( \half r(r+1) \eps h_j^{-1} \big)
  + (\bar \omega_{j+1} - \bar\omega_{j}) \big( \half r(r+1) \eps h_{j+1}^{-1} \big).
\end{align*}
If the mesh is uniform or if $r = 0$, then $\rho = 0$. In general,
under the mesh regularity assumption \eqref{eq:intro:kmesh_condition}
we obtain the residual estimate
\begin{equation}
  \label{eq:app:rho_residual}
  \|\rho\|_{\ell^\infty} \leq C(\kappa) \eps r.
\end{equation}
To estimate the error on the weights, we note that $M(\bar\omega -
\omega) = \rho$, and hence
\begin{displaymath}
  \| \bar\omega - \omega \|_{\ell^\infty} \leq
  \|M^{-1}\|_{\ell^\infty} \|\rho\|_{\ell^\infty} \leq
  \max(1,r)^{-1} C(\kappa) r \eps,
\end{displaymath}
that is, we obtain
\begin{equation}
  \label{eq:E:error_weights}
  \| \bar\omega - \omega\|_{\ell^\infty} \leq \cases{ C(\kappa) \eps\quad
    \text{if~} r \geq 1, \\
    0\quad \text{if~} r = 0. }
\end{equation}
This may seem an impossibly strong result at first glance, however, we
note that it is only true under the restriction that $r_k^\pm \equiv
r$ for all $k$. We expect that, in general, the cluster size will
influence the estimate to give an error of order $\mathcal{O}(\eps\max
r_k^\pm)$.

Upon noticing that the weights for the force-based summation rule
satisfy $\nu_j = \omega_j/\eps$, an estimate for the weights $\nu_j$
follows immediately.

\subsection{Proof of (\ref{eq:F:estimate_tildef})}
\label{app:int_err}
Let $I$ be the interpolation operator for the quasicontinuum mesh,
that is, $I : \X \rightarrow \X_h$ with $(I v_h)_{\ell_k} = v_{h,
  \ell_k}$, $k = -K+1, \dots, K$. Let $g \in \X$ be given by $g_\ell =
f_\ell \zeta_j(\eps \ell)$, then
\begin{displaymath}
  \sum_{\ell = -N+1}^N \eps g_\ell = \sum_{-N+1}^N \eps I g_\ell
  + \sum_{-N+1}^N \eps (g_\ell - Ig_\ell).
\end{displaymath}
In view of our assumption that $f_\ell = \bar f(\eps \ell),$ where
$\bar f \in \CC^2[-1,1]$, and the interpolation error estimate of
\cite[Thm A.4]{Ortner:2008a}, it follows that
\begin{displaymath}
  \sum_{\ell = -N+1}^N \eps g_\ell = \sum_{-N+1}^N \eps I g_\ell
  + \mathcal{O}(h_j^2 + h_{j+1}^2).
\end{displaymath}
Since, by definition, piecewise affine functions are summed exactly by
the cluster summation rule, it follows that
\begin{displaymath}
  \sum_{\ell = -N+1}^N \eps g_\ell = \sum_{k = -K+1}^K \nu_k
  \sum_{\ell \in \Cs_k}  \eps I g_\ell + \mathcal{O}(h_j^2 + h_{j+1}^2).
\end{displaymath}
Applying the same argument as above, we can deduce that
\begin{displaymath}
  \sum_{k = -K+1}^K \nu_k \sum_{\ell \in \Cs_k}  \eps I g_\ell
  = \sum_{k = -K+1}^K \nu_k \sum_{\ell \in \Cs_k}  \eps g_\ell
  + \mathcal{O}(h_j^2 + h_{j+1}^2),
\end{displaymath}
from which the desired result follows.

\end{document}